\documentclass[12pt]{amsart}

\usepackage{epsfig}
\usepackage{graphics}
\usepackage{latexsym}
\usepackage{psfrag}
\usepackage{amssymb}
\usepackage{amsmath}
\usepackage{amsthm} 

\newcommand{\cc}{{\mathcal{C}}}

\renewcommand{\setminus}{{\smallsetminus}}

\newcommand{\RR}{{\mathbb{R}}}

\newcommand{\NN}{{\mathbb{N}}}

\newcommand{\bigmid}{{~\mid~}}
\newcommand{\st}{{\bigmid}}

\newcommand{\nin}{{\notin}}

\newcommand{\closure}[1]{{\overline{#1}}}

\newcommand{\homeo}{{\medspace \cong \medspace}} 
 
\newcommand{\bdy}{{\partial}} 

\newcommand{\cross}{{\times}}

\newcommand{\interior}{{\rm{interior}}}

\theoremstyle{plain}
\newtheorem{theorem}{Theorem}[section]
\newtheorem{corollary}[theorem]{Corollary}
\newtheorem{lemma}[theorem]{Lemma}

\theoremstyle{definition}
\newtheorem*{define}{Definition}
\newtheorem{claim}[theorem]{Claim}
\newtheorem*{proofclaim}{Claim}

\newtheorem{remark}[theorem]{Remark}

\newtheorem*{question}{Question}

\newsavebox{\savepar}

\newcommand{\lab}[1]{{\ensuremath{\mathbb{#1}}}}
\newcommand{\piG}{{\pi_\Gamma}}

\begin{document}

\title{Surface Bundles versus Heegaard Splittings} 
\author{David Bachman}
\author{Saul Schleimer}
\address{Mathematics Department, Cal Poly State University, San Luis
  Obispo, CA 93407} 
\email{dbachman@calpoly.edu}
\address{Mathematics Department, University of Illinois at Chicago,
  Chicago, IL 60607}
\email{saul@math.uic.edu}

\date{\today}

\begin{abstract}
This paper studies Heegaard splittings of surface bundles via the
curve complex of the fibre.  The {\em translation distance} of the
monodromy is the smallest distance it moves any vertex of the curve
complex.  We prove that the translation distance is bounded above in
terms of the genus of any strongly irreducible Heegaard splitting.  As
a consequence, if a splitting surface has small genus compared to the
translation distance of the monodromy, the splitting is standard.
\end{abstract}
\maketitle

\section{Introduction}

The purpose of this paper is to show a direct relationship between the
action of a surface automorphism on the curve complex and the Heegaard
splittings of the associated surface bundle.  

We begin by restricting attention to automorphisms of closed
orientable surfaces with genus at least two.  Let $\varphi: F \to F$
be such an automorphism.  The {\em translation distance},
$d_\cc(\varphi)$, is the shortest distance $\varphi$ moves any vertex
in the curve complex.  (We defer precise definitions to
Section~\ref{Background}.)  As a bit of notation let $M(\varphi)$
denote the mapping torus of $\varphi$.  The following theorem gives a
link between translation distance and essential surfaces in
$M(\varphi)$.

\vspace{2mm}
\noindent
{\bf Theorem~\ref{IncompressibleBoundsTranslationDistance}.}  
{\em If $G \subset M(\varphi)$ is a connected, orientable,
  incompressible surface then either $G$ is isotopic to a fibre, $G$
  is homeomorphic to a torus, or $d_\cc(\varphi) \leq -\chi(G)$.}
\vspace{1mm}

An underlying theme in the study of Heegaard splittings is that, in
many instances, {\em strongly irreducible} splitting surfaces may take
the place of incompressible surfaces.  We prove:

\vspace{2mm}
\noindent
{\bf Theorem~\ref{StronglyIrreducibleBoundsTranslationDistance}.}
{\em If $H$ is a strongly irreducible Heegaard splitting of
  $M(\varphi)$ then $$d_\cc(\varphi) \leq -\chi(H).$$} 
\vspace{1mm}

This result gives new information about Heegaard splittings of
hyperbolic three-manifolds.  Previous work, by Y.~Moriah and
H.~Rubinstein~\cite{MoriahRubinstein97}, discussed the low genus
splittings of negatively-curved manifolds with very short geodesics.
Restricting attention to surface bundles and applying 
Theorems~\ref{StronglyIrreducibleBoundsTranslationDistance}
and~\ref{IncompressibleBoundsTranslationDistance} gives:

\vspace{2mm}
\noindent
{\bf Corollary~\ref{LowGenusImpliesStablized}.}  
{\em Any Heegaard splitting $H$ of the mapping torus $M(\varphi)$ with
  $-\chi(H) < d_\cc(\varphi)$ is a stabilization of the standard
  splitting.} 
\vspace{1mm}

This improves a result due, independently, to
H.~Rubinstein~\cite{Rubinstein02} and M.~Lackenby~\cite{Lackenby02}
(see our Corollary~\ref{HighPowersImpliesStandardSplittingUnique}).

The two theorems indicate an interesting connection between the
combinatorics of the curve complex and the topology of
three-manifolds.  This is in accordance with other work.  For example,
Y.~Minsky {\em et.al.} have used the curve complex to prove the
Ending Lamination Conjecture.  A major step is using a path in the
curve complex to give a model of the geometry of a hyperbolic
three-manifold.

Another example of this connection is found
in~\cite{CullerJacoRubinstein82} and~\cite{FloydHatcher82}.  These
papers study surface bundles where the fibre is a once punctured
torus.  Here the Farey graph takes the place of the curve complex.
An analysis of $\varphi$-invariant lines in the Farey graph
allows a complete classification of incompressible surfaces in such
bundles.

As we shall see in the proofs of
Theorems~\ref{StronglyIrreducibleBoundsTranslationDistance}
and~\ref{IncompressibleBoundsTranslationDistance}, essential surfaces
and strongly irreducible splittings in the mapping torus $M(\varphi)$
yield $\varphi$-invariant lines in the curve complex of the fibre.  It
is intriguing to speculate upon axioms for such lines which would,
perhaps, lead to classification results for essential surfaces or
strongly irreducible splittings.  This would be a significant step in
the over-all goal of understanding the topology of surface bundles.

At the heart of our proof of
Theorem~\ref{StronglyIrreducibleBoundsTranslationDistance} lies the
idea of a ``graphic'', due to Rubinstein and
Scharlemann~\cite{RubinsteinScharlemann96}. The graphic is obtained,
as in D.~Cooper and M.~Scharlemann's paper~\cite{CooperScharlemann99},
by comparing the bundle structure with a given height function and
applying Cerf theory.  As in their work, our situation requires a
delicate analysis of behavior at the vertices of the graphic.

The rest of the paper is organized as follows: basic definitions
regarding Heegaard splittings, surface bundles, and the curve complex
are found in Section~\ref{Background}.  With this background we
restate the main theorem and corollaries in Section~\ref{MainTheorem}.
Of main importance is the nature of simple closed curve intersections
between a fibre of the bundle and the Heegaard splitting under
discussion. This is covered in Section~\ref{AnalyzingIntersections},
in addition to a preliminary sketch of the proof of
Theorem~\ref{StronglyIrreducibleBoundsTranslationDistance}.  Our
version of the Rubinstein-Scharlemann graphic is discussed in
Section~\ref{Graphic}.  Concluding the paper
Section~\ref{ProvingMainTheorem} proves
Theorem~\ref{StronglyIrreducibleBoundsTranslationDistance} and poses a
few open questions.  We also discuss the possibility of strengthening
the inequality given in the main theorem.  For example, using a larger
``moduli space'' such as the pants complex or Teichmuller space does
not work.

We thank I.~Agol and M.~Culler for several helpful conversations.  In
particular Agol has shown us a construction, using techniques from his
paper~\cite{Agol02}, of surface bundles over the circle containing
strongly irreducible splittings of high genus. 

\section{Background material}
\label{Background}

This section presents the definitions used in this paper. A more
complete reference for the curve complex may be found
in~\cite{MasurMinsky99} while the paper~\cite{Scharlemann02} is an
excellent survey on Heegaard splittings.

\subsection{The curve complex}

Fix a closed orientable surface $F$ with genus $g(F) \geq 2$. 
If $\alpha \subset F$ is an essential simple closed curve then let
$[\alpha]$ be the isotopy class of $\alpha$.

\begin{define}
The set $\{[\alpha_0], \ldots, [\alpha_k]\}$ determines a {\em
$k$-simplex} if for all $i \neq j$ the isotopy classes $[\alpha_i]$,
$[\alpha_j]$ are distinct and there are $\alpha'_i \in [\alpha_i]$,
$\alpha'_j \in [\alpha_j]$ with $\alpha'_i \cap \alpha'_j =
\emptyset$.
\end{define}

\begin{define}
The {\em curve complex} of $F$ is the simplicial complex $\cc(F)$
given by the union of all simplices, as above.
\end{define}

We will restrict our attention to the zero and one-skeleta, $\cc^0(F)
\subset \cc^1(F)$.  Giving each edge length one the graph $\cc^1(F)$
becomes a metric space.  Let $d_\cc(\alpha, \beta)$ be the distance
between the vertices $[\alpha], [\beta] \in \cc^0(F)$.  When it can
cause no confusion we will not distinguish between an essential simple
closed curve and its isotopy class.

\begin{define}
Suppose that $\varphi$ is a homeomorphism of $F$.  The {\em
translation distance} of $\varphi$ is
$$d_\cc(\varphi) = \min \{ d_\cc(\alpha, \varphi(\alpha)) \st
\alpha \in \cc^0(F) \}.$$
\end{define}

\subsection{Heegaard splittings}

Fix a compact, orientable three-manifold $M$. Recall that a {\em
compression-body} is a three-manifold obtained as follows: Choose a
closed, orientable surface $H$ which is not a two-sphere.  Let $N = H
\cross I$.  Attach two-handles to the boundary component $H \cross \{
0 \}$.  Glue three-handles to all remaining boundary components which
are two-spheres.  We henceforth identify $H$ with the surface $H
\cross \{1\} \subset \bdy V$. 

If the resulting compression-body has only one boundary component then
it is a {\em handlebody}.

\begin{define}
A surface $H \subset M$ is a {\em Heegaard splitting} of $M$ if $H$ cuts
$M$ into a pair of compression-bodies $V$ and $W$.
\end{define}

\begin{define}
A properly embedded disk $D$ inside a compression-body $V$ is {\em
essential} if $\bdy D \subset H \subset \bdy V$ is essential.
\end{define}

Using this simple definition we arrive at an important notion, introduced
by A.~Casson and C.~Gordon~\cite{CassonGordon87}:

\begin{define}
A Heegaard splitting $H \subset M$ is {\em strongly irreducible} if all
pairs of essential disks $D \subset V$ and $E \subset W$ satisfy $\bdy
D \cap \bdy E \neq \emptyset$. 
\end{define}

If the splitting $H$ is not strongly irreducible then $H$ is {\em
  weakly reducible}.  M.~Scharlemann's~\cite{Scharlemann98} ``no
nesting'' lemma shows the strength of Casson and Gordon's definition:

\begin{lemma}
\label{NoNesting}
Suppose $H \subset M$ is a strongly irreducible splitting and $D
\subset M$ is an embedded disk with $\bdy D \subset H$ and with
$\interior(D)$ transverse to $H$.  Then there is a disk $D'$ properly
embedded in $V$ or $W$ with $\bdy D' = \bdy D$.
\end{lemma}

We take the following bit of terminology almost directly
from~\cite{ScharlemannThompson93}.  Let $M = F \cross I$ where $F$ is
a closed orientable surface.  Let $\alpha = \{ \rm{pt} \} \cross I$ be
a properly embedded arc.  Take $N$ a closed regular neighborhood of
$\bdy M \cup \alpha$ in $M$.  Let $H = \bdy N \setminus \bdy M$.  Then
$H$ is the {\em standard type 2} splitting of $M$.  The {\em standard
  type 1} splitting is isotopic to the surface $F \cross \{ 1/2 \}$.

Scharlemann and Thompson then prove:
\begin{theorem}
\label{SplittingsOfSurfaceCrossI}
Every Heegaard splitting of $M = F \cross I$ is a stabilization of
the standard type 1 or 2 splitting. 
\end{theorem}

Note that it is common to refer to stablizations of ``the'' standard
splitting as being standard themselves.

\subsection{Surface bundles}

Fix $F$, a closed, orientable surface with genus $g(F) \geq 2$.  Let
$\varphi : F \rightarrow F$ be a surface diffeomorphism which preserves
orientation.

\begin{define}
The {\em surface bundle} with monodromy $\varphi$ is the manifold
$$M(\varphi) = (F \cross [0, 2\pi])/\{(x, 2\pi) \equiv (\varphi(x), 0)\}.$$
\end{define}

Let $F(\theta)$ be the image of $F \cross \theta$. These surfaces are
{\em fibres} of the bundle $M(\varphi)$.  There is a natural smooth
map $\pi_F : M(\varphi) \rightarrow S^1$ defined by $\pi_F(x) =
\theta$ whenever $x \in F(\theta)$.  The map $\pi_F$ realizes
$M(\varphi)$ as a surface bundle over the circle.

We now define the {\em standard} Heegaard splitting of the surface
bundle $M(\varphi)$.  Pick $x, y \in F$ such that $x \neq y$ and
$\varphi(y) \neq x$. Fix $A$ and $B$ disjoint closures of regular
neighborhoods of $x \cross [0, \pi]$ and $y \cross [\pi, 2\pi]$
respectively.  Set
$$V = \closure{ (F \cross [\pi, 2\pi] \setminus B) \cup A}$$ and
$$W = \closure{ (F \cross [0, \pi] \setminus A) \cup B}.$$
Then $H = \bdy V = \bdy W$ is the standard Heegaard splitting of 
$M(\varphi)$.  Note that the genus of the standard splitting is $2g(F)
+ 1$.  Finally, the standard splitting is always weakly reducible
because the tubes $A \cap H$ and $B \cap H$ admit disjoint compressing
disks.

\section{Main theorem and corollaries}
\label{MainTheorem}

Let $F$ be a closed orientable surface with genus $g(F) \geq 2$.  Let
$\varphi: F \rightarrow F$ be an orientation-preserving
diffeomorphism.  The surface bundle $M(\varphi)$ is irreducible
and has minimal Heegaard genus two or larger.  Here then is a precise
statement of our main theorem:

\vspace{2mm}
\noindent
{\bf Theorem~\ref{StronglyIrreducibleBoundsTranslationDistance}.}
{\em If $H \subset M(\varphi)$ is a strongly irreducible Heegaard
splitting then the translation distance of $\varphi$ is at most the
negative Euler characteristic of $H$.  That is,
$$d_\cc(\varphi) \leq -\chi(H).$$} 
\vspace{1mm}

This is a deeper version of the following theorem
(from~\cite{Schleimer02a}): 

\begin{theorem}
\label{IncompressibleBoundsTranslationDistance}
If $G \subset M(\varphi)$ is a connected, orientable,
incompressible surface then either
\begin{itemize}
\item $G$ is isotopic to a fibre or
\item $G$ is a torus and $d_\cc(\varphi) \leq 1$ or
\item $d_\cc(\varphi) \leq -\chi(G)$.
\end{itemize}
\end{theorem}

For completeness we include a proof.

\begin{proof}[Proof of
Theorem~\ref{IncompressibleBoundsTranslationDistance}.] 

Suppose that $G$ is not isotopic to a fibre of the surface bundle.  If
$G$ is a torus then the map $\varphi$ is reducible.  It follows that
$d_\cc(\varphi) \leq 1$.  Assume, therefore, that $G$ is not a torus.

Briefly, the rest of the proof is as follows: Isotope $G$ into a
``good position'' and examine the transverse intersections of $G$ with
the fibres.  From these extract a sequence of curves which provide a
path in the curve complex of the fibre.  This path gives the desired
bound.  Here are the details.

Applying a theorem of W.~Thurston~\cite{Thurston86} isotope $G$ until
all non-transverse intersections with the fibres occur at a finite
number of saddle tangencies.  Furthermore, there is at most one
tangency between $G$ and any fibre $F(\theta)$.


It follows that {\em every} transverse curve of intersection between
$G$ and $F(\theta)$ is essential in both surfaces.  Any transverse
curve failing this would, perforce, be inessential in both.  But that
would lead directly to a center tangency between $G$ and some fibre.

Let $\{ \theta_i \}_{i = 0}^{n - 1}$ be the critical angles where $G$
fails to be transverse to $F(\theta_i)$.  Every critical angle gives a
saddle for $G$.  It follows that $n = - \chi(G) \geq 2$ as $G$ is
orientable and not a torus.  Pick regular angles $\{ \tau_i \}_{i =
  0}^{n - 1}$ where $\theta_{i - 1} < \tau_i < \theta_i$, with indices
taken mod $n$.  Apply a rotation to force $\tau_0 = 0 = 2\pi$.  Let
$\alpha_i$ be any curve component of $F(\tau_i) \cap G$ and recall
that $\alpha_i \subset F(\tau_i)$ is essential.

\begin{figure}
\psfrag{s}{\small$F(\tau_i)$}
\psfrag{t}{\small$F(\tau_{i+1})$}
\psfrag{a}{\small$\alpha_i$}
\psfrag{b}{\small$\alpha_{i+1}$}
\psfrag{c}{\small$\alpha_{i+1}''$}
$$\begin{array}{c}
\epsfig{file=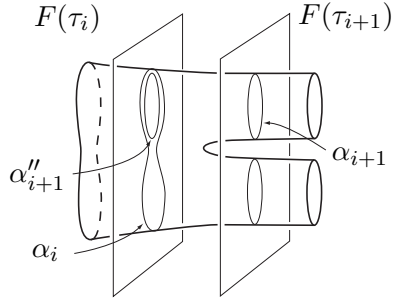, height = 3.8 cm}
\end{array}$$
\caption{Isotoping $\alpha_{i+1}$ back to the curve $\alpha_{i+1}''
  \subset F(\tau_i)$.}
\label{BackIsotopy}
\end{figure}

The curve $\alpha_{i+1}$ may be isotoped back through $F \cross
[\tau_i, \tau_{i+1}]$ to lie on the fibre $F(\tau_i)$, disjointly from
$\alpha_i$.  See Figure~\ref{BackIsotopy}.  In this fashion produce a
sequence of curves $\{ \alpha_i' \}_{i = 0}^n$ (indices {\em not}
taken mod $n$) in $F(\tau_0) = F(0)$ such that
\begin{itemize}
\item 
$\alpha_i'$ is isotopic to $\alpha_i$ through $F \cross [0,
\tau_{i}]$ for $i \in \{1, 2, \ldots, n - 1 \}$ while $\alpha_0' =
\alpha_0$ and $\alpha_n'$ is obtained by isotoping $\alpha_0$
off of $F(2\pi)$ back through $F \cross [0, 2\pi]$,

\item
$\alpha_i' \cap \alpha_{i+1}' = \emptyset$ for 
$i = 0, 1, \ldots, n - 1$, and

\item
$\varphi(\alpha_n')$ is isotopic to $\alpha_0'$.
\end{itemize}

Thus
$$d_\cc(\varphi) \leq d_\cc(\alpha_n', \varphi(\alpha_n')) 
                   \leq n = - \chi(G)$$
and the theorem follows. 
\end{proof}

Theorems~\ref{StronglyIrreducibleBoundsTranslationDistance}
and~\ref{IncompressibleBoundsTranslationDistance} have direct
consequences:

\begin{corollary}
\label{LowGenusImpliesStablized}
Any Heegaard splitting $H \subset M(\varphi)$ satisfying $-\chi(H) <
d_\cc(\varphi)$ is a stabilization of the standard splitting.
\end{corollary}

\begin{remark}
\label{HighDistanceImpliesStandard}
In particular, when the translation distance of $\varphi$ is bigger
than $4g(F)$, the standard splitting is the unique minimal genus
splitting, up to isotopy.
\end{remark}

\begin{corollary}
\label{HighPowersImpliesStandardSplittingUnique}
Suppose that $\varphi$ is pseudo-Anosov. If $n \in \NN$ is
sufficiently large then the standard Heegaard splitting of
$M(\varphi^n)$ is the unique splitting of minimal genus, up to
isotopy.
\end{corollary}

This follows from our
Corollary~\ref{LowGenusImpliesStablized} and
Lemma~4.6 of~\cite{MasurMinsky99}.

\begin{remark}
H.~Rubinstein~\cite{Rubinstein02} and M.~Lackenby~\cite{Lackenby02}
have independently and concurrently obtained
Corollary~\ref{HighPowersImpliesStandardSplittingUnique} using
techniques from minimal surface theory.  Note also that a closed
hyperbolic surface bundle $M(\varphi)$ is covered by the bundle
$M(\varphi^n)$.  Thus
Corollary~\ref{HighPowersImpliesStandardSplittingUnique} may be
considered as weak evidence for a ``yes'' answer to M. Boileau's
question, stated as Problem~3.88 in~\cite{Kirby97}.
\end{remark}

\begin{proof}[Proof of Corollary~\ref{LowGenusImpliesStablized}] 
Let $H \subset M(\varphi)$ be a Heegaard splitting with $-\chi(H) <
d_\cc(\varphi)$.  As the genus of the fibre $F$ is two or greater, the
genus of $H$ is at least two.  Thus the translation distance
$d_\cc(\varphi) > 2$ and $\varphi$ is irreducible.  By
Theorem~\ref{StronglyIrreducibleBoundsTranslationDistance} the
splitting $H$ cannot be strongly irreducible.  Thus $H$ is weakly
reducible. 

If $H$ is reducible then $H$ is stabilized as $M$ is
irreducible~\cite{Waldhausen68}.  In this case destabilize $H$ and
apply induction.  Suppose instead that $H$ is irreducible.
Following~\cite{CassonGordon87} there are disjoint disk systems ${\bf
  D} \subset V$ and ${\bf E} \subset W$ with the following property:
compressing $H$ along ${\bf D} \cup {\bf E}$ and deleting all
resulting two-sphere components yields a nonempty incompressible
surface $G$.  This surface need not be connected.  Also note that
$-\chi(G) < -\chi(H)$.

If a component of $G$ were a torus then either the genus of $F$ was
one or $\varphi$ was reducible.  Both possibilities yield
contradiction.  If a component of $G$ is not isotopic to the fibre
then apply Theorem~\ref{IncompressibleBoundsTranslationDistance} to
find $d_\cc(\varphi) \leq -\chi(G) < -\chi(H)$.  Again, this gives a
contradiction.

The only possibility remaining is that $G$ is isotopic to a collection
of fibres.  Note that the surface $G$ was obtained via compressing a
separating surface.  Thus $G$ itself must be separating.  Hence $G$ is
the union of an even number of fibres.  Thus the genus of $H$ is at
least that of the standard splitting.  (Note that even this much of
the proof establishes, when $d_\cc(\varphi) > 4g(F)$, that the standard
splitting has minimal genus.)

Let $N \homeo F \cross I$ be the closure of some component of $M
\setminus G$.  Let ${\bf A}$ be the compressing disks of ${\bf D}
\cup {\bf E}$ which meet $N$.  Finally, let $H' \subset N$ be the
surface obtained by compressing $H$ along all disks of $({\bf D}
\cup {\bf E}) \setminus {\bf A}$.

Then $H'$ is a Heegaard splitting for $N$ such that both boundary
components of $N$ lie on the same side of $H'$.  By
Theorem~\ref{SplittingsOfSurfaceCrossI} it follows that $H'$ is a
stabilization of standard type 2 splitting of $N \homeo F \cross I$.
The identical argument applies to all other components of $M(\varphi)
\setminus G$.  We conclude that $H$ is obtained by {\em amalgamation}
of these splittings (see~\cite{Schultens93}).  It follows that $H$
itself is a stabilization of the standard splitting of $M(\varphi)$
and this completes the proof.
\end{proof}

\section{Analyzing intersections}
\label{AnalyzingIntersections}

Briefly, the proof of
Theorem~\ref{StronglyIrreducibleBoundsTranslationDistance}  is as
follows: Isotope the splitting surface $H$ into a ``good position''
and examine the transverse intersections of $H$ with the fibres.  From
these extract a sequence of curves which provide a path in the curve
complex of the fibre.  This path gives the desired bound.  As might be
expected the details are more delicate than in the proof of
Theorem~\ref{IncompressibleBoundsTranslationDistance}.  In particular a
replacement for Thurston's theorem is needed. 

The rest of this section discusses the nature of intersections between
a Heegaard splitting surface and the fibres of a surface bundle over
the circle.

Fix a surface automorphism $\varphi: F \to F$.  Recall that $\pi_F :
M(\varphi) \to S^1$ is the associated map realizing $M(\varphi)$ as a
bundle over the circle.  Let $H \subset M(\varphi)$ be a Heegaard
splitting.  Pick a fibre $F(\theta) = \pi_F^{-1}(\theta)$ which meets
$H$ transversely.

\begin{define}
A simple closed curve component $\alpha \subset F(\theta) \cap H$ is
{\em non-compressing} if $\alpha$ is either essential in both
surfaces or inessential in both.
\end{define}


Note that, as $F(\theta)$ is incompressible, no curve of intersection
may be essential in $F(\theta)$ and also inessential in the splitting
surface.

\begin{define}
A simple closed curve component $\alpha \subset F(\theta) \cap H$ is
{\em mutually essential} if $\alpha$ is essential in both surfaces.
The curve $\alpha$ is {\em mutually inessential} if $\alpha$ is
inessential in both surfaces.
\end{define}

Wiggle $H$ slightly so that $\pi_F|H$ is a Morse function.  Let $p \in
H$ be a saddle critical point of $\pi_F|H$.  Let $c = \pi_F(p)$ and
let $P$ be the component of $H \cap \pi_F ^{-1}[c - \epsilon, c +
\epsilon]$ containing $p$.

\begin{define}
If every component of $\bdy P$ is mutually essential we call $p$ an
{\em essential} saddle.  If every component of $\bdy P$ is
non-compressing but at least one is inessential in $H$ then we call
$p$ an {\em inessential} saddle.
\end{define}

\begin{lemma}
\label{InterestingIntersection}
Fix $H \subset M(\varphi)$ so that $\pi_F|H$ is Morse.  If $F'$ is
isotopic to the fibre $F(\theta)$, the surface $H'$ is isotopic to
$H$, and $F'$ is transverse to $H'$ then $F' \cap H'$ contains at
least one curve which is not a mutually inessential curve.
\end{lemma}

\begin{proof}
Suppose that $F' \cap H'$ meet only in mutually inessential curves.
Then, by an innermost disk argument, there is a further isotopy
making them disjoint.  This cannot be, as handlebodies do not
contain closed embedded incompressible surfaces.
\end{proof}

\begin{lemma}
\label{LackOfMotion}
Fix $H \subset M(\varphi)$ so that $\pi_F|H$ is Morse.  Suppose angles
$\theta_- < \theta_+$ are given such that:
\begin{itemize}
\item 
The splitting $H$ intersects $F(\theta_\pm)$ transversely.

\item 
For every angle $\theta \in [\theta_-, \theta_+]$ all simple closed
curve components of $F(\theta) \cap H$ are non-compressing.

\item
Every saddle of $\pi_F|H$ in $F \cross [\theta_-, \theta_+]$ is inessential.
\end{itemize}

Then there is a curve of $F(\theta_-) \cap H$, essential in
$F(\theta_-)$, which is isotopic, through $F \cross [\theta_-,
\theta_+]$, to a curve of $F(\theta_+) \cap H$.
\end{lemma}

\begin{proof}
Let $\{\theta_i\}$ be the critical angles of $\pi_F|H$ which lie in
$[\theta_-, \theta_+]$.  Choose $r_i$ slightly greater than the $\theta_i$
and let $R = \{r_i\} \cup \{\theta_- + \epsilon\}$.  Refer to $R$ as
the set of {\em regular angles}.

For every $r \in R$ surger $H$ along every curve of $F(r) \cap H$
which bounds a disk in $F(r)$, innermost first.  Let $H'$ be the
intersection of the surgered surface with $F \cross [\theta_-,
\theta_+]$.  Note that $\pi_F|H'$ has exactly two new critical points
for every surgery curve.  See Figure~\ref{TrivialSurgery}.

\begin{figure}
\psfrag{H}{\small$H$}
\psfrag{h}{\small$H'$}
\psfrag{F}{\small$F(r)$}
$$\begin{array}{c}
\epsfig{file=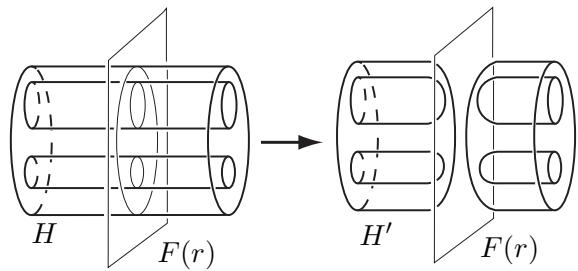, height = 3.5 cm}
\end{array}$$
\caption{Constructing $H'$ from $H$.}
\label{TrivialSurgery}
\end{figure}

\begin{proofclaim}
The surface $H'$ is a union of spheres, disks, and annuli. Every
annulus component has boundary which is essential in one of the
$F(\theta_\pm)$.
\end{proofclaim}

\begin{proof}
Fix attention on a component $H'' \subset H'$.  Begin by examining the
critical points of $\pi_F|H''$ and drawing conclusions about the Euler
characteristic of $H''$.

If $\pi_F|H''$ has no critical points then $H''$ is a horizontal
annulus.  In this case the boundary of $H''$ must be essential in its
fibre, by the construction of $H'$.  Note that this kind of situation
is the desired conclusion of the lemma at hand.

If $\pi_F|H''$ has more critical points of even index than odd then
$H''$ is a disk or sphere.

Now, suppose $p \in H''$ is a critical point of saddle type and let
$\pi_F(p) = \theta$.  Let $P$ be the component of $H \cap
\pi_F^{-1}[\theta - \epsilon, \theta + \epsilon]$ meeting $p$.

As $p$ is not an essential saddle at least one boundary component of
$P$ is inessential in its fibre.  If all three are inessential then
$H''$ is a two-sphere.  If exactly two components of $\bdy P$ are
inessential in their fibres then one of the two surfaces $F(\theta \pm
\epsilon)$ is compressible in $F \cross [\theta_-, \theta_+]$, a
contradiction.  Assume, therefore, that exactly one component of $\bdy
P$ is inessential in the containing fibre.  Call this curve $\alpha$.

Assume that $\pi_F(\alpha) = \theta - \epsilon$; that is, the
inessential curve $\alpha$ lies to the left of the saddle point. (The
other case is handled similarly.)  Let $r \in R$ be the regular value
appearing just before the critical value $\theta = \pi_F(p)$.  That
is, $r < \theta - \epsilon < \theta$.  Note that $\pi_F|H$ has no
critical values between $r$ and $\theta$.  See
Figure~\ref{AlphaSurgery}.

Deduce that, in the surface $H$, there is a horizontal annulus
isotopic to $\alpha \cross [r, \theta - \epsilon]$.  Thus the surface
$H''$ has a center tangency with the fibre just before the angle $r$.
This gives a disk capping off the curve $\alpha$.  Again, see
Figure~\ref{AlphaSurgery}.  It follows that every saddle in $H''$ is
paired with at least one critical point of even index.  Thus, if not a
sphere or disk, $H''$ is an annulus.

\begin{figure}
\psfrag{p}{\small$p$}
\psfrag{a}{\small$\alpha$}
\psfrag{c}{\small$F(\theta)$}
\psfrag{e}{\small$F(\theta-\epsilon)$}
\psfrag{r}{\small$F(r)$}
$$\begin{array}{c}
\epsfig{file=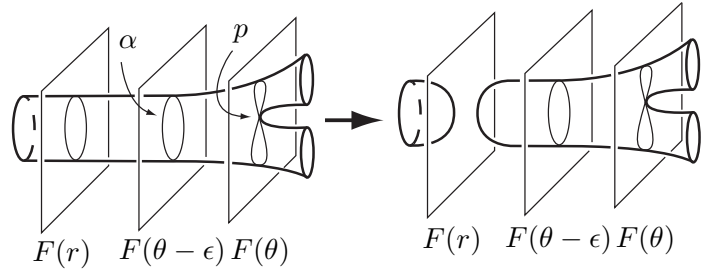, height = 3.5 cm}
\end{array}$$
\caption{A surgery just before $\alpha$}
\label{AlphaSurgery}
\end{figure}

Finally, if $H''$ is an annulus then $H''$ has two boundary
components.  By construction, each of these is an essential curve
in one of the surfaces $F(\theta_\pm)$.  This finishes the claim.
\end{proof}

Now to complete the proof of the lemma.  Suppose that no component of
$H'$ meets both boundary components of $F \cross [\theta_-,
\theta_+]$.  Thus, by the claim, every component of $H'$ meeting
$F(\theta_-)$ is boundary parallel in $F \cross [\theta_-, \theta_+]$.
Isotope $F(\theta_-)$ across these boundary-parallelisms to obtain a
surface $F'$ which intersects the splitting surface $H$ only in
mutually inessential curves.  This contradicts
Lemma~\ref{InterestingIntersection}.

So there is a component $H'' \subset H'$ which meets both
$F(\theta_-)$ and $F(\theta_+)$.  By the claim above this $H''$ must
be isotopic to a horizontal annulus with boundary essential in the
containing fibres.  The lemma is proved.
\end{proof}

\section{The graphic and its labellings}
\label{Graphic}

To begin, this section discusses height functions subordinate to a
Heegaard splitting.  We then compare one of these to a bundle
structure to obtain a ``graphic'' in the sense
of~\cite{RubinsteinScharlemann96}.  This technique is similar that
of~\cite{CooperScharlemann99}. We also refer to~\cite{KobayashiSaeki00} as
an informative paper on this topic. 

Fix attention now on a Heegaard splitting $H$ of a closed, orientable,
three-manifold $M$.  Recall that $H$ cuts $M$ into a pair of
handlebodies $V$ and $W$.

\subsection{Height functions}
Choose a diffeomorphism between the handlebody $V$ and a regular
neighborhood of a connected, finite, polygonal graph $\Theta \subset
\RR^3$.  For simplicity assume that every vertex of $\Theta$ has
valence two or three.  Let $\Theta_V$ be the image of $\Theta$ inside
of $V$.  Any such $\Theta_V$ is a {\em spine} of $V$.  

\begin{define}
A smooth map $\pi_H: M \rightarrow I$ is a {\em height function} with
respect to the splitting $H$ if
\begin{itemize}
\item 
the {\em level} $H(t) = \pi_H^{-1}(t)$ is isotopic to $H$ for all $t
\in (0,1)$,  
\item 
the graphs $\Theta_V = \pi_H^{-1}(0)$, $\Theta_W = \pi_H^{-1}(1)$
are spines for $V$ and $W$,
\item 
there is a map $h: H \cross I \rightarrow M$ (a {\em sweep-out}) such that 
\begin{itemize}
\item $h | H \cross (0,1)$ is a diffeomorphism, 
\item $\pi_H \circ h$ is projection onto the second factor,
\item for small $\epsilon$ the image of $h | H \cross [0, \epsilon]$ gives
the structure of a mapping cylinder to $V(\epsilon) =
\pi_H^{-1}[0,t]$, for some deformation retraction $\bdy V(\epsilon)
\to \Theta_V$,  
\item the previous condition also holds for the handlebody 
$W(1 - \epsilon) = \pi_h^{-1}[1-\epsilon,1]$.
\end{itemize}
\end{itemize}
\end{define}

The last two conditions on the sweep-out ensure that if an embedded
surface $F \subset (M, H)$ meets $\Theta_V$ (or $\Theta_W$)
transversely then $F \cap V(\epsilon)$ (or $F \cap W(1-\epsilon)$) is
a collection of properly embedded disks.  

\subsection{The graphic}
Let $H \subset M(\varphi)$ be a Heegaard splitting of the given
surface bundle.  Choose a bundle map, $\pi_F$, and a height function,
$\pi_H$, as above and insist that the two functions are generic with
respect to each other.  Define the map $\piG : M(\varphi) \to S^1 \cross
I$ by $\piG(x) = (\pi_F(x), \pi_H(x))$.  As a bit of terminology, we
sometimes call $\piG(F(\theta))$ a {\em vertical arc} and call $\piG(H(t))$
a {\em horizontal circle}.  Also, let $\Gamma(\theta, t) = \piG^{-1}(\theta,
t) = F(\theta) \cap H(t)$.  

The {\em graphic} of the map $\piG = (\pi_F, \pi_H)$ is the set
$$\Lambda = \overline{ \{ (\theta, t) \in S^1 \cross (0,1) 
\st \mbox{$F(\theta)$ is not transverse to $H(t)$} \}}$$
where the closure is taken in $S^1 \cross I$.  As in the papers cited
above: $\Lambda$ is a graph with smooth edges meeting the boundary of
the annulus transversely.

\begin{define}
The {\em regions} of the graphic are the open cells of $S^1 \cross
(0,1) \setminus \Lambda$.  
\end{define}

If $(\theta, t)$ is a point of a region then $\Gamma(\theta, t) =
F(\theta) \cap H(t)$ is a collection of simple closed curves embedded
in $M(\varphi)$.

\begin{remark}
\label{RegionsAreBoring}
Suppose that $(\theta, t), (\theta', t')$ are two points in the same
region.  Then the {\em combinatorics} of $\Gamma(\theta, t)$ and
$\Gamma(\theta', t')$ are identical, as follows:


Suppose that $\alpha$ is a short horizontal or vertical arc embedded in the
interior of a region $R$.  Suppose also that $(\theta, t), (\theta',
t')$ are the endpoints of $\alpha$.  Let $\gamma$ be a component of
$\Gamma(\theta, t)$.  Then there is an ambient isotopy taking $\gamma$
to $\gamma' \subset \Gamma(\theta', t')$ supported in a neighborhood
of an annulus of $F(\theta)$ (if $\alpha$ is vertical) or an annulus
of $H(t)$ (if $\alpha$ is horizontal).  Furthermore this ambient
isotopy may be chosen to take $F(\theta)$ to $F(\theta')$ and $H(t)$
to $H(t')$

By crawling along horizontal and vertical arcs any pair of points in
$R$ may be joined.  Thus most properties of $\Gamma(\theta, t)$ depend
only on the region containing $(\theta, t)$.
\end{remark}

As above, the {\em edges} are the one-dimensional strata of $\Lambda$.

\begin{remark}
\label{EdgesAreInteresting}
When $(\theta, t)$ lies on an edge there is one component,
$\Sigma(\theta, t)$, of $\Gamma(\theta, t)$ which is not a simple
closed curve.  This $\Sigma$ is the {\em singular component}.  The
name of the point $(\theta, t)$ is omitted when clear from context.

There are two kinds of edge: those representing a center tangency
between a fibre and a level surface and those representing a saddle
tangency.  Crossing an edge of the graphic from region $R$ to $R'$
causes the combinatorics of the curves to change.  If the edge
represents a center tangency then a single curve of $R$ disappears
(appears).  
It follows that this curve is mutually inessential; in
other words it bounds a disk in both the fibre and the level.  In this
situation $\Sigma$ is a single point.

If the edge represents a saddle then two curves in $R$ touch as the
edge is crossed and become a single curve in $R'$ (or the reverse).
When $(\theta, t)$ lies on such an edge the singular component
$\Sigma$ is a four valent graph with one vertex, embedded in both
the fibre $F(\theta)$ and the level $H(t)$. 
\end{remark}

The {\em vertices} are the zero-dimensional strata of the graphic
$\Lambda$.  

\begin{figure}
\psfrag{H(t)}{\small$\theta$}
\psfrag{F(theta)}{\small$t$}
$$\begin{array}{cc}
\epsfig{file=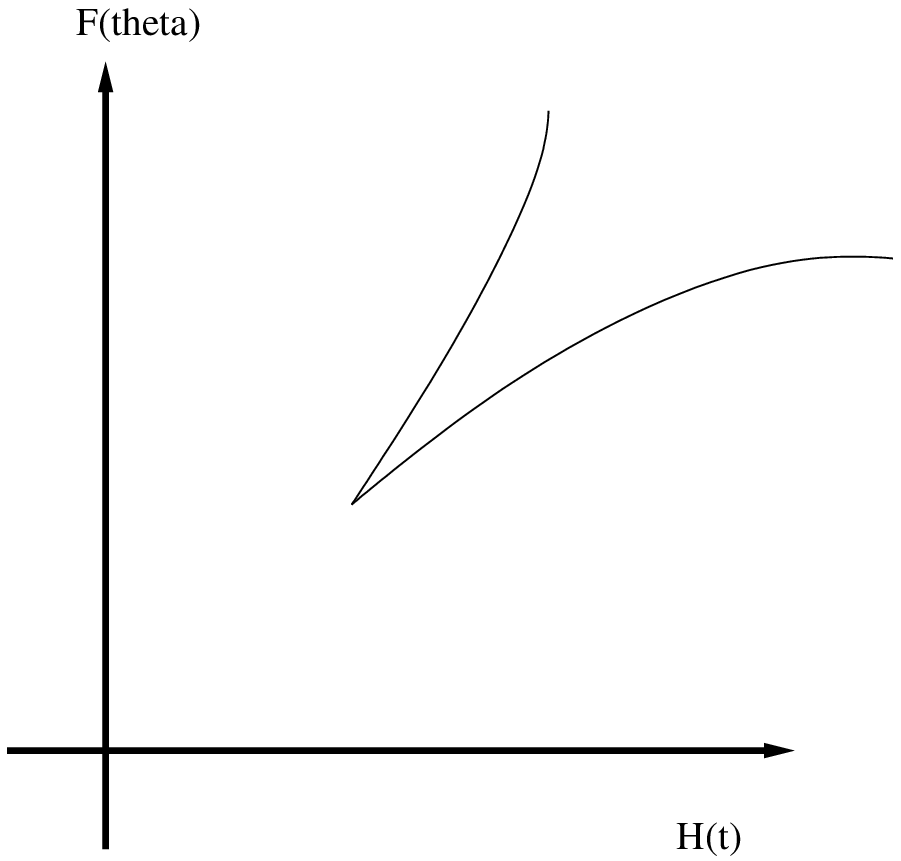, height = 3.5 cm} &
\epsfig{file=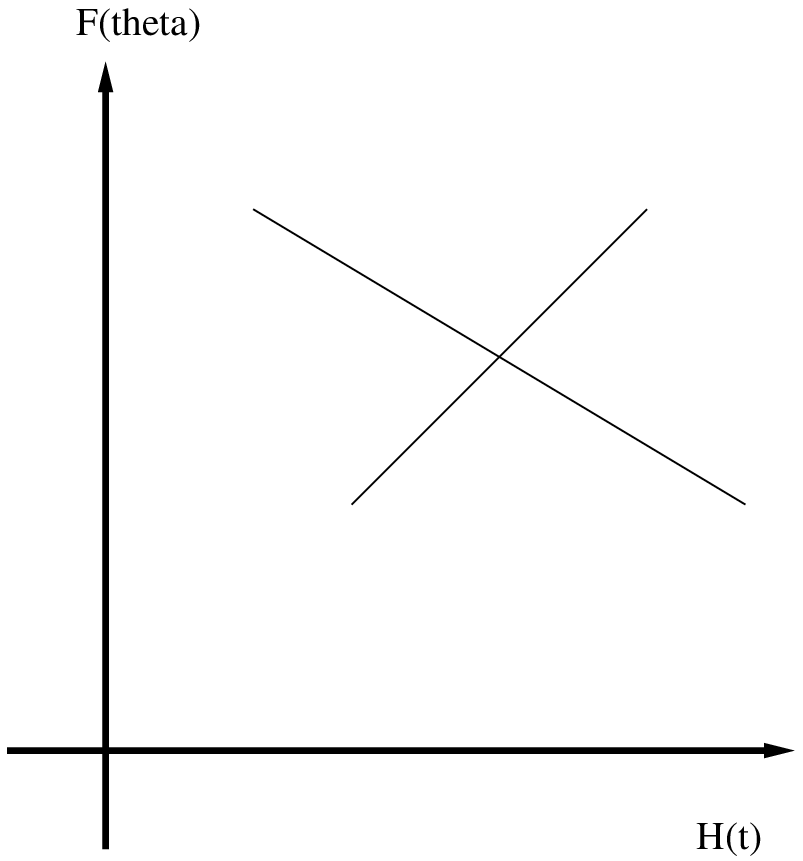, height = 3.5 cm} \\
\end{array}$$
\caption{A birth-death vertex and a crossing vertex}
\label{VerticesInGraphic}
\end{figure}

\begin{remark}
\label{VerticesAreInteresting}
There are several possibilities for a vertex $v = (\theta, t)$:
\begin{enumerate}
\item
All vertices of valence one occur at height $0$ or $1$.
\item 
A vertex with valence two in $\Lambda$ is a {\em birth-death}
vertex. Both edges lie in the same quadrant with respect to the
vertex.  (As in~\cite{KobayashiSaeki00}.)  See
Figure~\ref{VerticesInGraphic}.
\item 
A vertex of valence four is a {\em crossing} vertex. Here the
four edges lie in distinct quadrants and the tangent directions of
opposite edges agree.  The edges cut a small neighborhood of the
vertex into four regions.  Again see Figure~\ref{VerticesInGraphic}.
There are two further subcases:
\begin{itemize}
\item 
If the singular component $\Sigma \subset \Gamma(\theta, t)$ has
two components then $v$ is a {\em disjoint crossing} vertex.
\item
If $\Sigma$ has only one component then $v$ has {\em entangled saddles}.  
In this case $\Sigma$ is a four valent graph with exactly two
vertices.
\end{itemize}
\end{enumerate}
\end{remark}

Finally, general position implies that for every angle $\theta$ the
vertical arc $\piG(F(\theta)) = \{\theta\} \cross I$ meets the graphic
$\Lambda$ at most once nontransversely, either at a vertex or at a
tangency with the interior of an edge.  The same holds for horizontal
circles $\piG(H(t)) = S^1 \cross \{t\}$.

\subsection{Labellings}

The following labellings will play an important role in the proof of
Theorem~\ref{StronglyIrreducibleBoundsTranslationDistance}.  Recall
that the level surface $H(t)$ cuts $M(\varphi)$ into a pair of handlebodies
$V(t)$ and $W(t)$.

\begin{define}
For each $t \in (0, 1)$ the level $H(t)$ is labelled \lab{V} (or
\lab{W}) if there is a value $\theta \in [0, 2\pi]$ and a simple
closed curve component $\gamma \subset \Gamma(\theta, t) = F(\theta)
\cap H(t)$ which bounds an essential disk in $V(t)$ (resp.  $W(t)$).
\end{define}

A slightly finer labelling is also required.

\begin{define}
A region $R$ is labelled with a \lab{V} (or \lab{W}) if there is a
$(\theta, t) \in R$ and a simple closed curve component $\gamma
\subset \Gamma(\theta, t)$ such that $\gamma$ bounds an essential
disk in $V(t)$ ($W(t)$).
\end{define}

As we shall see in the proof below each level and region receives at
most one label (or none at all) and these labels reveal a great deal
of information about relative positions of the fibre and level
surface.

\section{Proof of the main theorem}
\label{ProvingMainTheorem}

This section, after reiterating the statement, proves our main
theorem.

\begin{theorem}
\label{StronglyIrreducibleBoundsTranslationDistance}
If $H \subset M(\varphi)$ is a strongly irreducible Heegaard splitting
then the translation distance of $\varphi$ is at most the negative Euler
characteristic of $H$.  That is, $$d_\cc(\varphi) \leq -\chi(H).$$
\end{theorem}

Pick a bundle map $\pi_F : M(\varphi) \to S^1$ and a height function
$\pi_H : M(\varphi) \to I$ which respect the given fibre and
splitting.  Recall that $\Theta_V = \pi_H^{-1}(0)$ and
$\Theta_W = \pi_H^{-1}(1)$ are spines for the handlebodies $V$
and $W$, respectively, while the surfaces $H(t) = \pi_H^{-1}(t)$
impose a product structure on $M(\varphi) \setminus (\Theta_V \cup
\Theta_W)$.

As above obtain a graphic $\Lambda$ in the annulus $S^1 \cross I$.
We begin by obtaining a few fairly standard facts about the labellings
defined in the previous section.

\subsection{Analyzing the labelling}
\label{FindingGoodLevel}

\begin{claim}
\label{LabelsAreNonempty}
For all sufficiently small positive values $\epsilon$ the level
$H(\epsilon)$ is labelled $\lab{V}$ while $H(1 - \epsilon)$ is
labelled $\lab{W}$. 
\end{claim}

This follows directly from the construction of the height function and
genericity.  

\begin{claim}
\label{LabelsAreDisjoint}
No level or region is labelled with both a \lab{V} and a \lab{W}.
Also, if a region $R$ is labelled then every level $H(t)$, such that
$\piG(H(t)) \cap R \neq \emptyset$, receives the same label.  Finally,
if a level $H(t)$ is labelled then some region meeting $\piG(H(t))$
receives the same label.
\end{claim}

This follows from strong irreducibility, the fact that the
curves $\Gamma(\theta, t) = F(\theta) \cap H(t)$ form a singular
foliation of $H(t)$, and Remark~\ref{RegionsAreBoring}.

\begin{remark}
\label{LevelWitnessesAreStable}
Suppose that $H(t)$ is labelled and $\gamma \subset \Gamma(\theta, t)$
is a witness of this fact.  Then there is a $\gamma' \subset
\Gamma(\theta', t)$ which is also a witness, for all $\theta'$
sufficiently close to $\theta$.
\end{remark}

\begin{remark}
\label{NoNasties}
Suppose $H(t)$ is labelled with a \lab{V} (or \lab{W}). By
incompressibility of $F(\theta)$ the given curve $\gamma$ bounds a
disk $D \subset F(\theta)$.  By the ``no nesting'' lemma,
Lemma~\ref{NoNesting}, $D$ is isotopic rel $\gamma$ to a disk
properly embedded in $V(t)$ (resp.  $W(t)$).
\end{remark}

From Remark~\ref{EdgesAreInteresting} deduce:

\begin{claim}
\label{CenterTangencyEdge}
If two regions are both adjacent to an edge representing a center
tangency, then both regions have the same label.
\end{claim}

The next claim is not strictly required for the proof of
Theorem~\ref{StronglyIrreducibleBoundsTranslationDistance}.  We
include it both to simplify our proof of the theorem and to shed light
on the general situation.

\begin{claim}
\label{LabelConsistency}
Suppose that $H(t')$ is labelled \lab{V}. Suppose $0 < t < t'$. Then
$H(t)$ is labelled \lab{V} as well.  Identically, if $t' < t < 1$ and
$H(t')$ is labelled \lab{W}, then $H(t)$ is labelled \lab{W}.
\end{claim}

\begin{proof}
Consider the case where $H(t')$ is labelled \lab{V}.  By hypothesis we
are given an angle $\theta'$ and a curve $\gamma'' \subset
\Gamma(\theta', t')$ bounding an essential disk in $V(t')$.  We may
choose an angle $\theta$ close to $\theta'$ so that firstly, by
Remark~\ref{LevelWitnessesAreStable}, there is a curve $\gamma'
\subset \Gamma(\theta, t')$ bounding an essential disk in $V(t)$ and
secondly, by general position, the point $(\theta, t)$ lies in a
region of the graphic.  That is, $(\theta, t) \nin \Lambda$.  

By Remark~\ref{NoNasties} the curve $\gamma'$ bounds a disk $D'
\subset F(\theta)$.  Furthermore, choosing a different $\gamma'
\subset \Gamma(\theta, t')$ if necessary, assume that $\gamma'$ is
``innermost.''  That is, all curves of $\interior(D') \cap H(t')$ are
inessential in $H(t')$.

Let $\Gamma(D', t) = D' \cap H(t) \subset \Gamma(\theta, t)$.  Suppose
that all components of $\Gamma(D', t)$ are inessential in $H(t)$.
Isotope $D'$ rel $\gamma'$ off of $H(t)$ to a disk $E$ with $E \cap
V(t) = \emptyset$.  Let $N = \closure{V(t') \setminus V(t)}$ and note
that $N \homeo H \cross I$ while $\bdy N = H(t) \cup H(t')$.  Thus we
may further isotope $E$ rel $\gamma'$ out of $\interior(N)$.  This
pushes $E$ to a disk $E' \subset W(t')$.  Thus $\gamma'$ either bounds
an essential disk in $W(t')$ or is trivial in $H(t')$.  The first
implies that $H$ is reducible while the second directly contradicts
the hypothesis of the claim.

Thus there is a curve, $\gamma \subset \Gamma(D', t)$, which is
essential in $H(t)$ and is the innermost such in $D'$.  Let $D \subset
D'$ be the disk which $\gamma$ bounds.  Recall that all curves of
$\interior(D) \cap \bdy N$ are inessential in $\bdy N$.  So there is
an isotopy of $D$ rel $\gamma$ to a disk lying inside of $M(\varphi)
\setminus (H(t') \cup H(t))$.

Now, if this disk lies in $N = \closure{V(t') \setminus V(t)}$ then
$\gamma$ could not be essential in $H(t)$.  It follows that $\gamma$
bounds a disk in $V(t)$.  Thus $H(t)$ is labelled with a \lab{V}, as
claimed.
\end{proof}

As a bit of notation set 
$L(V) = \{ t \in (0, 1) \st H(t) \mbox{ is labelled } \lab{V} \}$ and  
$L(W) = \{ t \in (0, 1) \st H(t) \mbox{ is labelled } \lab{W} \}$.
Define $t_V = \sup L(V)$ and $t_W = \inf L(W)$

\begin{claim}
\label{DisjointOpenIntervals}
The sets $L(V)$ and $L(W) \subset (0, 1)$ are nonempty, disjoint,
connected, and open.  Also, $t_V \leq t_W$.
\end{claim}

\begin{proof}
The first sentence follows from Claims~\ref{LabelsAreNonempty},
\ref{LabelsAreDisjoint}, \ref{LabelConsistency}, and the fact that if
$\gamma \subset \Gamma(\theta, t)$ is a simple closed curve then
there are $\gamma' \subset \Gamma(\theta, t - \epsilon)$ and $\gamma''
\subset \Gamma(\theta, t + \epsilon)$ with combinatorics identical to
$\gamma$.  The second follows from the fact that $L(V)$, $L(W)$ are
disjoint open intervals and that $\inf L(V) = 0$ by
Claim~\ref{LabelsAreNonempty}.  
\end{proof}

It follows that if $t_V \leq t_0 \leq t_W$ then $H(t_0)$ is
unlabelled.  

\subsection{Analyzing the unlabelled level}

The unlabelled level $H(t_0)$ found above will serve as a replacement
for Thurston's theorem used in the proof of
Theorem~\ref{IncompressibleBoundsTranslationDistance}.  There are two
cases to consider: $t_V < t_W$ or $t_V = t_W$.  Each is dealt with, in
turn, below.

\subsection{Unlabelled interval}

Suppose that $t_V < t_W$.  Pick a level $H(t_0)$ which avoids the
vertices of the graphic, which is not tangent to any edge of the
graphic, and has $t_V < t_0 < t_W$.  It immediately follows that for
every $\theta$ and for every nonsingular $\gamma \subset
\Gamma(\theta, t_0)$ the curve $\gamma$ is non-compressing --- either
essential in both $F(\theta)$ and $H(t_0)$ or inessential in both.  In
the terminology of
Rubinstein-Scharlemann~\cite{RubinsteinScharlemann96} the level
$H(t_0)$ is {\em compression-free} with respect to every fibre.

We now come to the heart of this case: Let $\{ \theta_i \}_{i =
  0}^{n-1}$ be the critical angles of $\pi_F|H(t_0)$ corresponding to
essential saddles.  (As defined in
Section~\ref{AnalyzingIntersections}, the saddle point $p$ is {\em
  essential} when all boundary components of the associated pair of
pants are mutually essential curves.)  Choose the indexing so that
$\theta_i < \theta_{i + 1}$ are adjacent.  

Suppose first that $n = 0$.  Rotate the $S^1$ coordinate, if
necessary, so that $F(0)$ meets $H(t_0)$ transversely.  Cut
$M(\varphi)$ along $F(0)$.  Note that $H(t_0)$, inside of $F \cross
[0, 2\pi]$, satisfies the hypotheses of Lemma~\ref{LackOfMotion}.
Thus there is an mutually essential curve of $F(0) \cap H(t_0)$
isotopic, through $F \cross [0, 2\pi]$, to a curve of $F(2\pi)$.  Thus
$d_\cc(\varphi) \leq 1 < 2 \leq -\chi(H)$, as desired.

Suppose now that $n > 0$.  Choose $\epsilon$ sufficiently small and
positive.  So, for each critical angle, there is an essential saddle
in $H(t_0) \cap (F \cross [\theta_i - \epsilon, \theta_i +
\epsilon])$.  The pair of pants given by this saddle contributes $-1$
to the Euler characteristic of $H(t_0)$.  Thus $n \leq -\chi(H)$.

By Lemma~\ref{LackOfMotion}, there are also mutually essential curves
$\alpha_i \subset \Gamma(\theta_i + \epsilon, t_0)$ and $\alpha_i'
\subset \Gamma(\theta_{i+1} - \epsilon, t_0)$ such that $\alpha_i'$ is
isotopic to $\alpha_i$ through $F \cross [\theta_i + \epsilon,
\theta_{i + 1} - \epsilon]$.  Also $\alpha_{i}$ may be isotoped back
through $F \cross [\theta_i - \epsilon, \theta_i + \epsilon]$ to lie
in $F(\theta_i - \epsilon)$, disjoint from $\alpha_{i-1}'$.  (This is
shown by Figure~\ref{BackIsotopy}, although the labels will be
different.)

As in the proof of
Theorem~\ref{IncompressibleBoundsTranslationDistance} the $\alpha_i$
give a path of length $n$ in $\cc^1(F(0))$ starting at $\alpha_0$ and
ending at $\varphi(\alpha_0)$.  This implies $d_\cc(\varphi) \leq
-\chi(H)$, as desired.

\subsection{Only a vertex}
The more difficult situation occurs when $t_0 = t_V = t_W$.  By our
general position assumption the horizontal circle $\piG(H(t_0))$ meets
the graphic $\Lambda$ at most once nontransversely.

\begin{claim}
\label{ExistsAVertex}
There are regions $R$ and $R'$ labelled $\lab{V}$ and $\lab{W}$,
respectively, such that the closures $\closure{R}$ and $\closure{R'}$
both meet $\piG(H(t_0))$.  Also, neither $R$ nor $R'$ meet
$\piG(H(t_0))$.  The horizontal circle $\piG(H(t_0))$ meets a crossing
vertex $(\theta_0, t_0)$ of the graphic.
\end{claim}

\begin{proof}
As $t_0 = t_V$ and by Claim~\ref{LabelsAreDisjoint} there is a
region $R$ such that $\closure{R} \cap \pi_\Gamma(H(t_0))$ is
nonempty and $R$ lies below $\pi_\Gamma(H(t_0))$ in the annulus $S^1
\cross I$.  Similarly there is a region $R'$ above the horizontal
circle $\pi_\Gamma(H(t_0))$

Now, if $\piG(H(t_0))$ is transverse to the edges of $\Lambda$ then
every region $R''$ with closure meeting $\piG(H(t_0))$ also has
interior meeting $\piG(H(t_0))$.  Then, as $H(t_0)$ is unlabelled, so
are $R$ and $R'$.  This is a contradiction.

Suppose that $\piG(H(t_0))$ is only tangent to an edge of the graphic.
Then every region, but one, whose closure meets the circle
$\piG(H(t_0))$ also meets $\piG(H(t_0))$ along its interior.  As
above, this gives a contradiction.  Thus $\piG(H(t_0))$ meets a
vertex at the point $(\theta_0, t_0)$.

Finally, to rule out the possibility that the vertex is a birth-death
vertex: When $\piG(F(\theta_0))$ and $\piG(H(t_0))$ meet in a
birth-death vertex there are only two edges of the graphic incident on
the vertex.  As in Figure~\ref{VerticesInGraphic} the slopes of the
two edges have the same sign and we may assume that, as the edges
leave the vertex, both edges head ``northeast''.  (The other three
cases are similar.)  Again, every region, but one, whose closure
meets the circle $\piG(H(t_0))$ also meets $\piG(H(t_0))$ along its
interior.  Again, this is a contradiction.
\end{proof}

Now focus attention on this vertex at $(\theta_0, t_0)$.  Let $\Sigma
\subset \Gamma(\theta_0, t_0)$ be the union of the singular
components.  Let $P$ be the components of $H(t_0) \cap (F \cross
[\theta_0 + \epsilon, \theta_0 - \epsilon])$ meeting $\Sigma$.  We
will call $P$ a {\em foliated regular neighborhood} of $\Sigma$, taken
in $H(t_0)$.  Similarly, let $Q$ be a foliated regular neighborhood of
$\Sigma$, taken in $F(\theta_0)$.  Finally, let $\bdy_\pm P = P \cap
F(\theta_0 \pm \epsilon)$ while $\bdy_\pm Q = Q \cap H(t_0 \pm
\epsilon)$.

\begin{claim}
The vertex at $(\theta_0, t_0)$ has entangled saddles.  Also, some
component $ \beta \subset \bdy_- Q$ bounds an essential disk in
$V(t_0 - \epsilon)$ (and some component $\delta \subset \bdy_+ Q$
bounds in $W(t_0 + \epsilon)$).
\end{claim}

\begin{proof}
As $\piG(F(\theta_0))$ and $\piG(H(t_0))$ meet in a crossing vertex
there are four regions adjacent with closure meeting the vertex.  Call
these ``north'', ``east'', ``south'', and ``west''.  See
Figure~\ref{VerticesInGraphic}. Again, our general position assumption
ensures that $\piG(H(t_0))$ meets the graphic $\Lambda$ at most once
nontransversely.  Thus $\piG(H(t_0))$ meets $\Lambda$ exactly once
nontransversely.  Thus all regions whose closure meets $\piG(H(t_0))$,
other than the north and south, also meet $\piG(H(t_0))$ along their
interior.  It follows from Claim~\ref{LabelsAreDisjoint} that all
these regions except the south (the region $R$) and north ($R'$) are
unlabelled.

Choose a curve $\beta \subset \Gamma(\theta_0, t_0 - \epsilon)$ which
bounds an essential disk in $V(t_0 - \epsilon)$.  Moving along a
straight arc from $(\theta_0, t_0 - \epsilon)$ to $(\theta_0 +
\epsilon, t_0)$ cannot induce ambient isotopy on $\beta$.  If it did
the east region would be labelled with a \lab{V}, an impossibility.

Thus, the southeast edge represents a saddle tangency.  Also, this
saddle meets a regular neighborhood of $\beta$, taken inside of
$F(\theta_0)$.  Symmetrically, the same holds for the southwest edge.
We conclude that both saddles lie inside the same component of
$\Gamma(\theta_0, t_0)$.  An identical argument locates $\delta \subset
\Gamma(\theta_0, t_0 + \epsilon)$.  It follows that $\Sigma$ is connected,
contains two saddle tangencies, and has $\beta$ and $\delta$ as
boundary components of $Q$, the vertical foliated neighborhood.
\end{proof}

Note that the curves $\beta$ and $\delta$ may be isotoped to curves
$\beta_0$ and $\delta_0$ lying inside of $P$, the regular neighborhood
of $\Sigma \subset H(t_0)$.  Since $\beta$ and $\delta$ bound disks in
$V(t_0 - \epsilon)$ and $W(t_0 + \epsilon)$ the curves $\beta_0$ and
$\delta_0$ bound disks in $V(t_0)$ and $W(t_0)$.

\begin{remark}
\label{PAndQDontFill}
Recall that the Heegaard splitting $H$ and the fibre $F$ both have
genus at least two.  The graph $\Sigma$ has only two vertices and four
edges.  Thus, by an Euler characteristic argument, there is an
essential simple closed curve in $H(t_0)$ disjoint from the subsurface
$P$ (and similarly for $Q \subset F(\theta_0)$).  So if $t_V = t_W$
then the Heegaard splitting $H$ satisfies Thompson's {\em disjoint
  curve property}, defined in~\cite{Thompson99}.
\end{remark}

\begin{figure}
\psfrag{P}{\small$P$}
\psfrag{S}{\small$\Sigma$}
\psfrag{a}{\small$\beta_0$}
\psfrag{b}{\small$\delta_0$}
\psfrag{a'}{\small$\beta'$}
\psfrag{b'}{\small$\delta'$}
\psfrag{a''}{\small$\beta''$}
\psfrag{b''}{\small$\delta''$}
$$\begin{array}{cc}
\epsfig{file=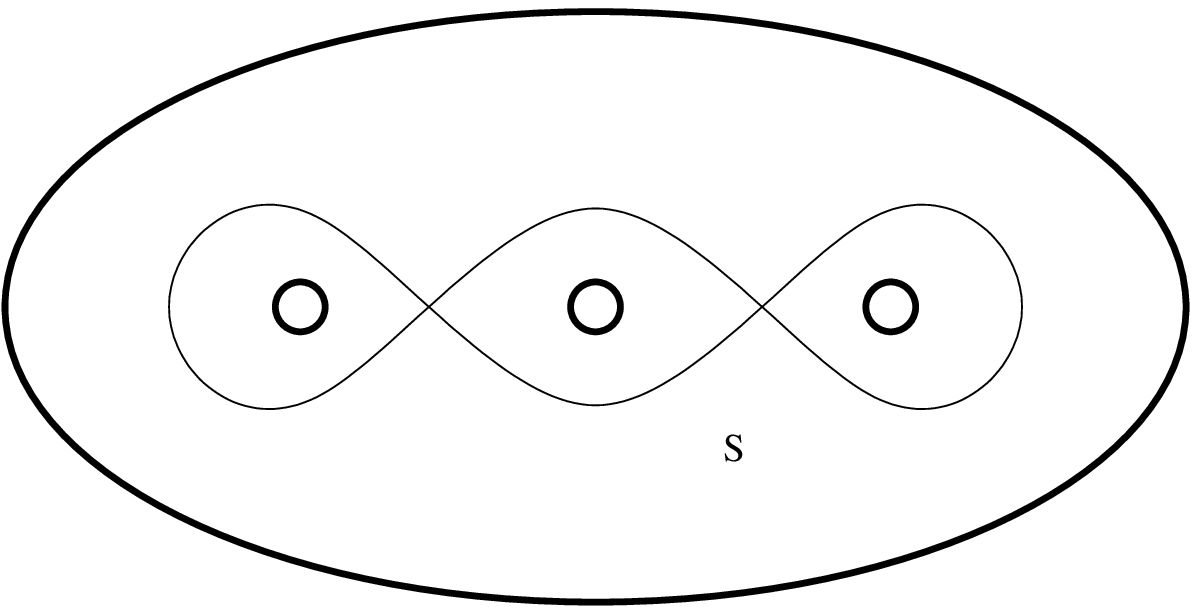, height = 2 cm} &
\epsfig{file=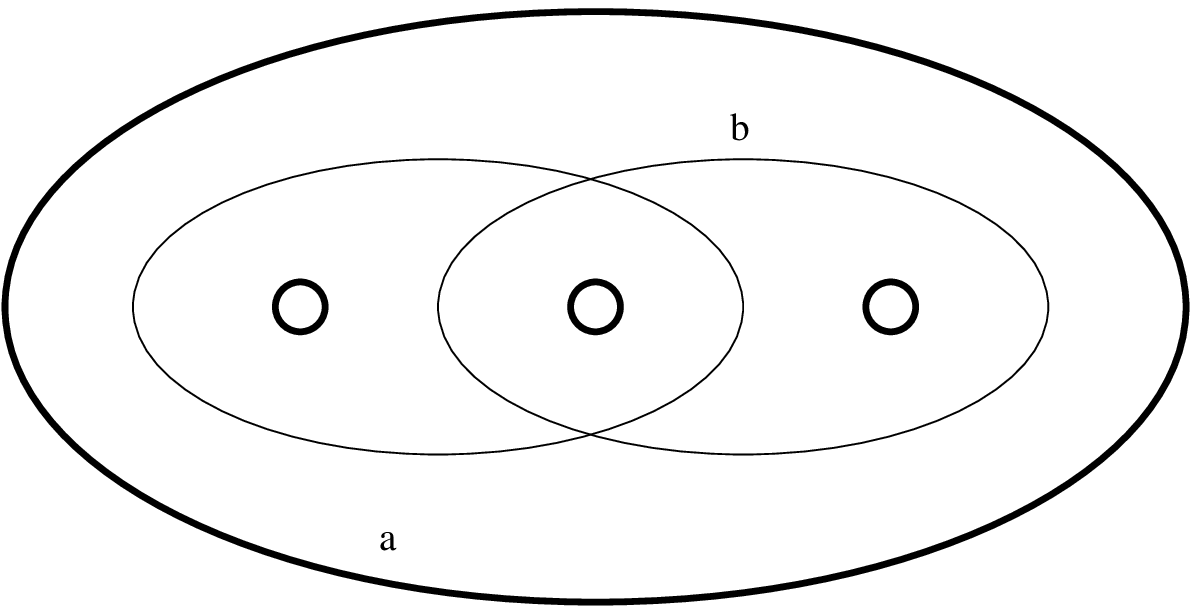, height = 2 cm} \\
\epsfig{file=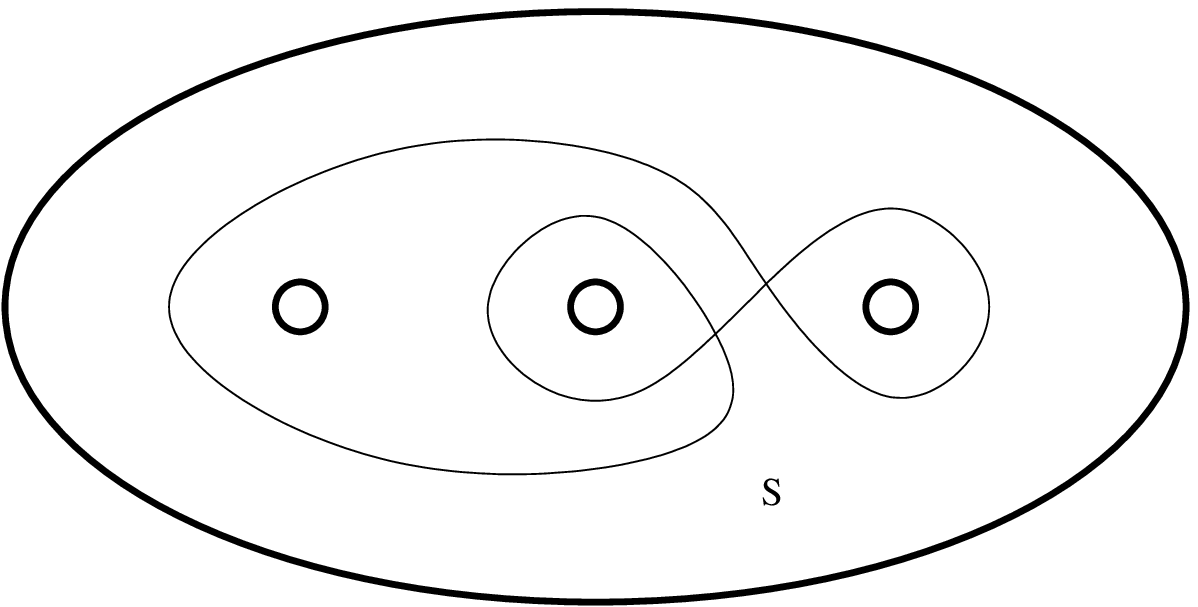, height = 2 cm} &
\epsfig{file=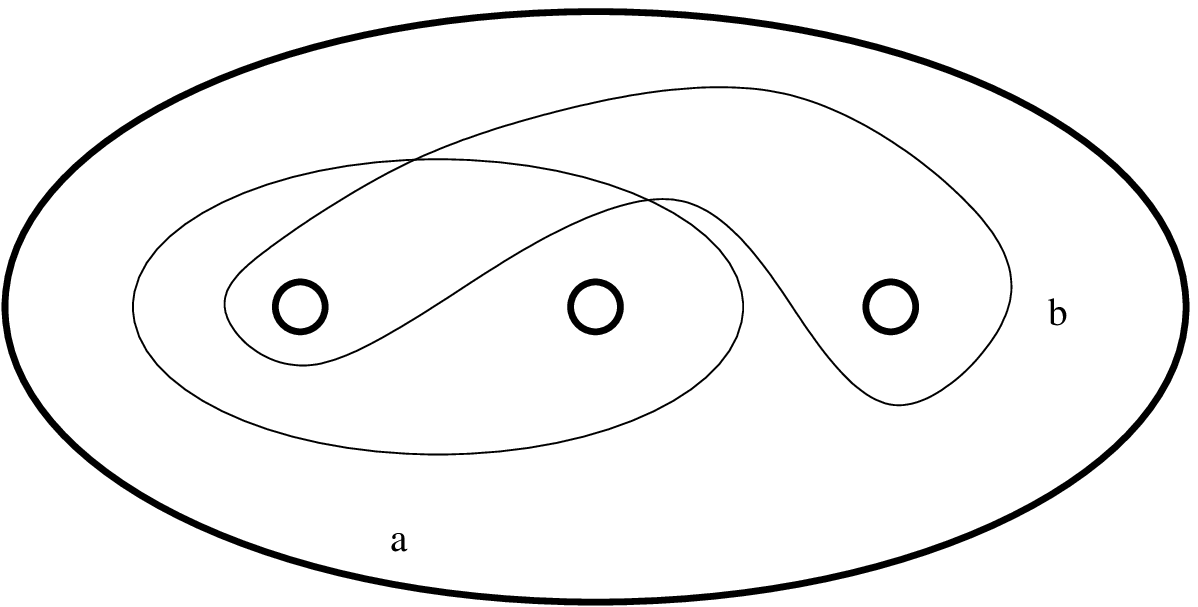, height = 2 cm} \\
\epsfig{file=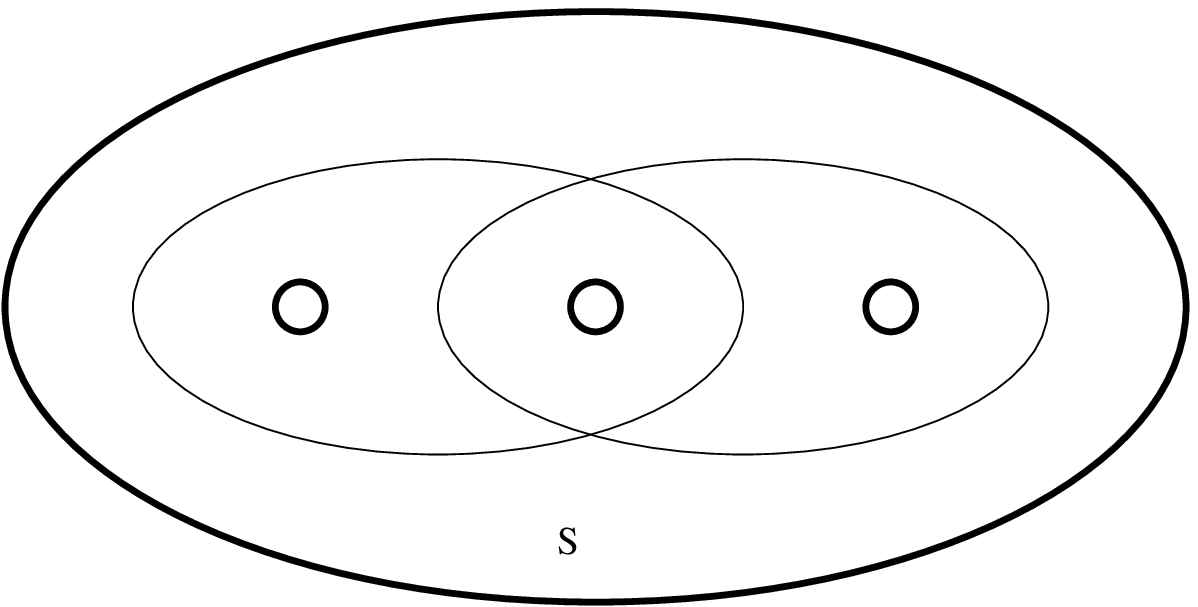, height = 2 cm} &
\epsfig{file=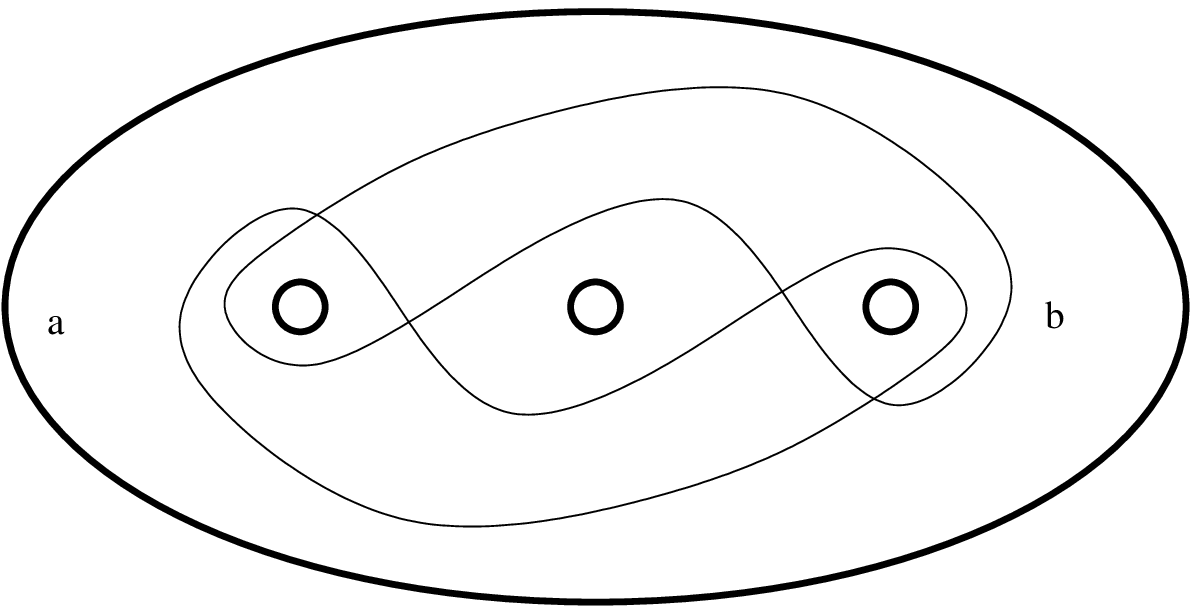, height = 2 cm} \\
\epsfig{file=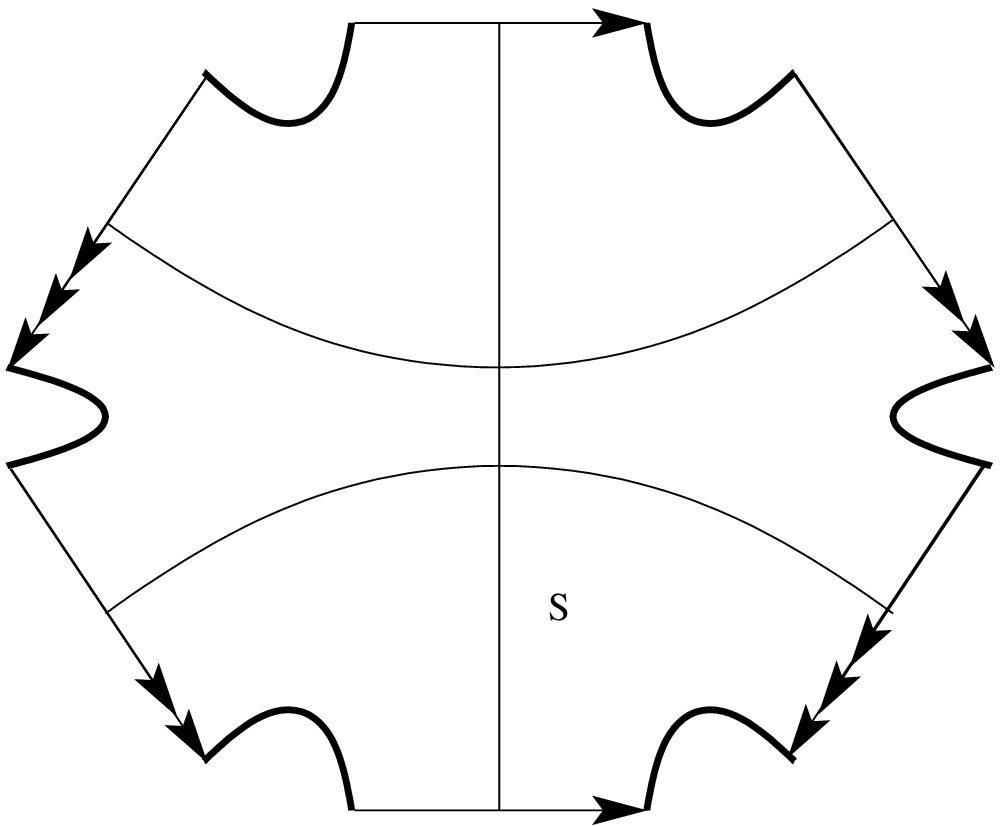, height = 3 cm} &
\epsfig{file=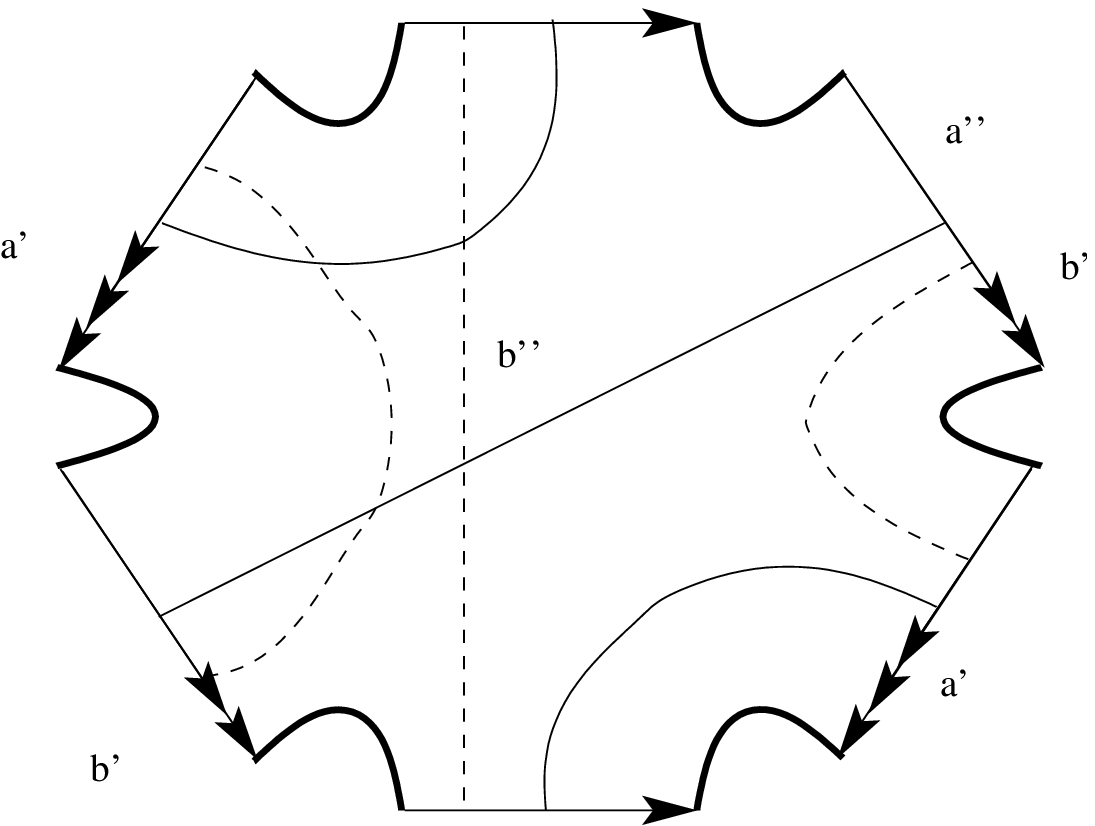, height = 3 cm} \\
\end{array}$$
\caption{Possibilities for the regular neighborhood $P$.}
\label{PNeighborhood}
\end{figure}

We now examine the properties of the foliated neighborhoods of
$\Sigma$; that is, $P \subset H(t_0)$ and $Q \subset F(\theta_0)$.
Taking into account symmetry, orientability of $H$, and the
co-orientation of the foliation itself there are four possibilities
for $P$, as shown in the first column of Figure~\ref{PNeighborhood}.
Note that there are four ways of resolving the two saddles of
$\Sigma$.  These, after performing a small isotopy, yield the curves
$\bdy_\pm P$ and $\bdy_\pm Q$.  Recall that $\beta \subset \bdy_- Q$
and $\delta \subset \bdy_+Q$ may be isotoped to $\beta_0, \delta_0
\subset P$, as shown in the second column of
Figure~\ref{PNeighborhood}.

In the fourth row, there are two possibilities for $\beta_0$; either
$\beta_0$ agrees with $\beta'$ or with $\beta''$ (the solid curves).
The curve $\delta_0$ is treated similarly (see the dashed curves).
Wherever $\beta_0$ and $\delta_0$ lie inside the twice-punctured
torus, we see that $\beta_0$ only meets $\delta_0$ once.  Thus the
splitting surface $H$ is stabilized, contradicting strong
irreducibility.  It follows that the foliated neighborhood $P$ is
homeomorphic to a four-times punctured sphere, with $\Sigma$ in
various positions, as shown.

\begin{claim}
Every component of $\bdy P$ is essential in $H(t_0)$.
\end{claim}

\begin{proof}
In each of the three cases if any component of $\bdy P$ bounds a disk
in $H(t_0)$ then a component of $\bdy P$ bounds a disk in
$\closure{H(t_0) \setminus P}$.  Isotope $\beta_0$ across this disk to
make $\beta_0$ disjoint from $\delta_0$.  This contradicts the strong
irreducibility of $H$.
\end{proof}

\begin{claim}
\label{PEssential}
The components of $\bdy P$ are mutually essential curves.
\end{claim}

\begin{proof}
All are essential in $H(t_0)$ by the above claim.  If one of the
curves bounds a disk in the containing fibre then, by
Lemma~\ref{NoNesting}, that curve bounds a disk in $V(t_0)$ or
$W(t_0)$.  In this case $H(t_0)$ was labelled with a $\lab{V}$ or a
$\lab{W}$.  This contradicts Claim~\ref{DisjointOpenIntervals} and
the hypothesis $t_0 = t_V = t_W$.
\end{proof}

Now to carry out an analysis similar to that for the case $t_V < t_W$.
Let $\{ \theta_i \}_{i = 1}^{n-1}$ be the angles corresponding to the
essential saddles of $\pi_F|H(t_0)$.  Let $(\theta_0, t_0)$ be the
vertex with entangled saddles meeting the horizontal circle
$\piG(H(t_0))$.  Choose the indexing so that $\theta_i < \theta_{i +
  1}$ are adjacent.  Choose $\epsilon$ sufficiently small and
positive.  So for each critical angle, $i > 0$, there is an essential
saddle in $H(t_0) \cap (F \cross [\theta_i - \epsilon, \theta_i +
\epsilon])$ which contributes $-1$ to the Euler characteristic of
$H(t_0)$.  As the curves of $\bdy P$ are mutually essential
(Claim~\ref{PEssential}) the four-punctured sphere $P$ contributes
$-2$ to the Euler characteristic of $H(t_0)$.  Altogether, the $n - 1$
essential saddles and $P$ contribute $n + 1$ to the negative Euler
characteristic of $H(t_0)$.  That is, $n + 1 \leq -\chi(H)$.

By Lemma~\ref{LackOfMotion}, there are also mutually essential curves
$\alpha_i \subset \Gamma(\theta_i + \epsilon, t_0)$ and $\alpha_i'
\subset \Gamma(\theta_{i+1} - \epsilon, t_0)$ such that $\alpha_i'$ is
isotopic to $\alpha_i$ through $F \cross [\theta_i + \epsilon,
\theta_{i + 1} - \epsilon]$.  Also, for $i > 0$, $\alpha_{i}$ may be
isotoped back through $F \cross [\theta_i - \epsilon, \theta_i +
\epsilon]$ to lie in $F(\theta_i - \epsilon)$, disjoint from
$\alpha_{i-1}'$.  For $i = 0$ this last may not hold.  Instead, after
isotoping $\alpha_0$ through $F \cross [\theta_0 - \epsilon, \theta_0
+ \epsilon]$ to obtain $\alpha_0''$, both $\alpha_0''$ and $\alpha_{n
  - 1}'$ may lie in a translate of $Q$, the vertical foliated
neighborhood of $\Sigma$.  As in Remark~\ref{PAndQDontFill} the
subsurface $Q$ does not fill $F(\theta_0 - \epsilon)$.  In any case
$d_\cc(\alpha_0'', \alpha_{n-1}) \leq 2$ when considered in the
curve complex of the fibre.

Thus, similar to the proof of
Theorem~\ref{IncompressibleBoundsTranslationDistance}, the $\alpha_i$
give a path of length $n + 1$ in the graph $\cc^1(F(0))$.  Therefore
$d_\cc(\varphi) \leq n + 1 \leq -\chi(H)$, as desired.  This deals
with the case where $t_V = t_W$ and proves the theorem.

We end here with a few remarks and comments:

\begin{question}
What can be said about the higher genus Heegaard splittings of
surface bundles with high translation distance?
\end{question}

\begin{question}
Suppose that $M$ has two distinct surface bundle structures. What can be
learned from the graphic induced on $S^1 \cross S^1$?
\end{question}

\begin{question}
It seems likely that the techniques of Section~\ref{FindingGoodLevel}
are sufficiently soft to allow a taut foliation to replace the surface
bundle structure.  See also Question~9.5 of D.~Calegari's problem list
on foliations~\cite{Calegari02}.
\end{question}

\begin{remark}
There is no inequality, as in Theorem~\ref{MainTheorem}, between the
genus of a strongly irreducible splitting of $M(\varphi)$ and the
stretch factor of the pseudo-Anosov automorphism $\varphi$.  Here is
the required construction: Fix $H \subset M(\varphi)$ a genus two,
strongly irreducible splitting.  (For example, let $M(\varphi)$ be the
longitudinal filling on D.~Rolfsen's $6_2$ knot~\cite{Rolfsen90}.  As
the knot is tunnel number one let $H$ be the resulting genus two
Heegaard splitting.)

Isotope $H$ until all curves of $F(0) \cap H$ are mutually essential.
Let $\gamma \subset \Gamma = F(0) \cap H$ be one component.  Let
$N(\gamma)$ be a regular neighborhood of $\gamma$.  Let $\mu, \lambda
\subset \bdy N(\gamma)$ be a meridian, longitude pair with $\mu$
bounding a disk in $N(\gamma)$ while $\lambda$ is isotopic to $\bdy
N(\gamma) \cap H$ (which, in turn, is isotopic to $\bdy N(\gamma) \cap
F$).  Then $M_n$, the $1/n$ Dehn surgery on $\closure{M \setminus
  N(\gamma)}$, is still a surface bundle, with monodromy $\varphi_n$,
say.  The stretch factor of $\varphi_n$ grows linearly with $n$
(see~\cite{LongMorton86}) while the Heegaard genus of $M_n$ remains
equal to two.  Thus the minimal genus splitting remains strongly
irreducible.  This completes the construction.
\end{remark}

\begin{remark}
It has been asked whether the main theorem of this paper can be
improved to refer to translation distance in the pants complex.
(See~\cite{Brock01}.)  It is straight-forward to provide candidate
counterexamples in genus two, somewhat similar to the above.  A subtle
argument, shown to us by D.~Canary and Y.~Minsky, then proves that the
volumes increase without bound.  We plan on providing the details of
this construction in a future paper.
\end{remark}


\bibliographystyle{plain}
\bibliography{bibfile}
\end{document}